\documentclass[reqno,centertags]{amsart}
\usepackage{amsmath}
\usepackage{amssymb}
\usepackage{amsfonts}

\newtheorem{theorem}{Theorem}
\theoremstyle{plain}

\newtheorem{lemma}{Lemma}

\newtheorem{problem}{Problem}

\numberwithin{equation}{section}

\input{tcilatex}

\begin{document}
\title[]{Ruled surfaces of finite type with respect to the second fundamental form}
\author{Hassan Al-Zoubi }
\address{Department of Mathematics, Al-Zaytoonah University of Jordan, P.O.
Box 130, Amman, Jordan 11733}
\email{dr.hassanz@zuj.edu.jo}
\author{Amer Dababneh}
\address{Department of Mathematics, Al-Zaytoonah University of Jordan}
\email{dababneh.amer@zuj.edu.jo}
\author{Waseem Mashaleh}
\address{Department of Mathematics, Al-Zaytoonah University of Jordan}
\email{w.almashaleh@zuj.edu.jo}
%\author{Mutaz Al-Sabbagh}
%\address{Department of Basic Sciences and Humanities, Imam Abdulrahman bin Faisal University}
%\email{malsbbagh@iau.edu.sa}
\author{Nancy Ramahi}
\address{Department of Basic Sciences, Al-Zaytoonah University of Jordan}
\email{N.Alramahi@zuj.edu.jo}
\date{}

\begin{abstract}
In this article, we consider surfaces in the 3-dimensional Euclidean space $\mathbb{E}^{3}$ without parabolic points which are of finite $II$-type, that is, they are of finite type, in the sense of B.-Y. Chen, corresponding to the second fundamental form. We study an important family of surfaces, namely, ruled surfaces in $\mathbb{E}^{3}$. We show that ruled surfaces are of infinite $II$-type.
\end{abstract}

\maketitle

\section{Introduction}
Euclidean immersions of finite type were introduced by B.-Y. Chen about four decades ago and it has been a topic of active research by many differential geometers since then. Many results concerning this subject can be found in \cite{B3}.
A submanifold $M^{m}$ is said to be of finite type corresponding to the first fundamental form $I$, if each component of the position vector field $\boldsymbol{x}$ of $M^{m}$ can be expressed as a finite sum of eigenfunctions of the Laplacian $\Delta ^{I}$, that is,
\begin{equation}
\boldsymbol{x}=\boldsymbol{c}+\sum_{i=1}^{k}\boldsymbol{x}_{i},\ \ \ \ \
\label{0}
\end{equation}%
where $\Delta ^{I}\boldsymbol{x}_{i}=\beta _{i}\boldsymbol{x}_{i},i=1,...,k$, $\boldsymbol{c}$ is a constant vector and $\beta _{1},\beta _{2},...,\beta _{k}$ are eigenvalues of $\Delta ^{I}$. Moreover, if there are exactly $k$ nonconstant eigenvectors $\boldsymbol{x}_{1},...,\boldsymbol{x}_{k}$ appearing in (\ref{0}) which all belong to different eigenvalues $\beta _{1},\beta _{2},...,\beta _{k}$, then $M^{m}$ is said to be of $I$-type $k$. However, if $\beta _{i}= 0$ for some $i = 1, ..., k$, then $M^{m}$ is said to be of null $I$-type $k$ , otherwise $M^{m}$ is said to be of infinite type.

%$\boldsymbol{x}_{0}$ is a fixed vector and $\boldsymbol{x}_{1},...,% \boldsymbol{x}_{k}$ are nonconstant maps such that If, in particular, all eigenvalues $\lambda _{1},\lambda _{2},...,\lambda _{k}$ are mutually distinct, then $S$ is said to be of $J$-type $k$, otherwise $S$ is said to be of infinite type. When $\lambda _{i}=0$ for some \emph{i} = 1,..., \emph{k%}, then $S$ is said to be of null $J$-type $k$.

In general when $M^{m}$ is of finite type $k$, it follows from (\ref{0}) that there exist a monic polynomial, say $R(x)\neq 0,$ such that $R(\Delta ^{I})(\boldsymbol{x}-\boldsymbol{x}_{0})=\mathbf{0}.$ Suppose that $R(x)=x^{k}+\sigma_{1}x^{k-1}+...+\sigma _{k-1}x+\sigma _{k},$ then coefficients $\sigma _{i}$ are given by

\begin{eqnarray}
\sigma _{1} &=&-(\lambda _{1}+\lambda _{2}+...+\lambda _{k}),  \notag \\
\sigma _{2} &=&(\lambda _{1}\lambda _{2}+\lambda _{1}\lambda_{3}+...+\lambda_{1}\lambda _{k}+\lambda _{2}\lambda _{3}+...+\lambda _{2}\lambda
_{k}+...+\lambda _{k-1}\lambda _{k}),  \notag \\
\sigma _{3} &=&-(\lambda _{1}\lambda _{2}\lambda _{3}+...+\lambda_{k-2}\lambda _{k-1}\lambda _{k}),  \notag \\
&&.............................................  \notag \\
\sigma _{k} &=&(-1)^{k}\lambda _{1}\lambda _{2}...\lambda _{k}.  \notag
\end{eqnarray}

Therefore the position vector $\boldsymbol{x}$ satisfies the following equation, (see \cite{B3})

\begin{equation}  \label{1}
(\Delta^{I})^{k}\boldsymbol{x}+\sigma_{1}(\Delta^{I})^{k-1}\boldsymbol{x}+...+\sigma_{k}( \boldsymbol{x}-\boldsymbol{x}_{0})=\boldsymbol{0}.
\end{equation}

The class of finite type submanifolds in an arbitrary dimensional Euclidean space is very large, meanwhile very little is known about surfaces of finite type in the Euclidean 3-space corresponding to the fundamental form $I$. Actually, so far, the only known surfaces of finite type in $\mathbb{E}^{3}$ are the minimal surfaces, the circular cylinders and the spheres. So in \cite{B2} B.-Y. Chen asked the following geometric quastion

%Let $\boldsymbol{x},\boldsymbol{H}$ be the position vector field and the mean curvature field of $M^{n}$ respectively.
%\begin{equation}
%\Delta ^{I}\boldsymbol{x}=-n\boldsymbol{H}.  \notag
%\end{equation}

%From this formula one can see that $M^{n}$ is a minimal submanifold if and only if all coordinate functions, restricted to $M^{n}$, are eigenfunctions of $\Delta ^{I}$ with eigenvalue $\lambda =0$. Moreover in \cite{R15} T. Takahashi showed that the submanifold $M^{n}$ for which $\Delta ^{I}% \boldsymbol{x}=\lambda \boldsymbol{x}$, i.e., for which all coordinate functions are eigenfunctions of $\Delta ^{I}$ with the same eigenvalue $%\lambda \in \mathbb{Re}$, are precisely either the minimal submanifold with eigenvalue $\lambda =0$ or the minimal submanifold of hyperspheres $S^{m-1}$ with eigenvalue $\lambda >0$. Although the class of finite type submanifolds in an arbitrary dimensional Euclidean spaces is very large, very little is known about surfaces of finite type in the Euclidean 3-space $E^{3}$. Actually, so far, the only known surfaces of finite type corresponding to the first fundamental form in the Euclidean 3-space are the minimal surfaces, the circular cylinders and the spheres. So in \cite{R5} B.-Y. Chen mentions the following problem

\begin{problem}
\label{(1)}Classify all surfaces of finite Chen $I$-type in $\mathbb{E}^{3}$.
\end{problem}

With the aim of getting an answer to this problem, important families of surfaces were studied by different authors by proving that finite type ruled surfaces \cite{B5}, finite type quadrics \cite{B6}, finite type tubes \cite{B4}, finite type cyclides of Dupin \cite{D1, D2}, finite type cones \cite{G1}, and finite type spiral surfaces \cite{B1} are the only known examples of surfaces in $\mathbb{E}^{3}$. However, for other classical families of surfaces, such as surfaces of revolution, translation surfaces as well as helicoidal surfaces, the classification of its finite type surfaces is not known yet.

In this area, S. Stamatakis and H. Al-Zoubi restored attention to this theme by introducing the notion of surfaces of finite type with respect to the second or third fundamental forms (see \cite{S1}). As an extension of the above problem, we raise the following two questions which seem to be very interesting:

\begin{problem}
\label{2}Classify all surfaces of finite $II$-type in the Euclidean 3-space.
\end{problem}

\begin{problem}
\label{3}Classify all surfaces of finite $III$-type in the Euclidean 3-space.
\end{problem}

Therefore, in order to give an answer to the second and third problem, it is worthwhile investigating the classification of surfaces in the Euclidean space $\mathbb{E}^{3}$ in terms of finite $J$-type, $(J=II, III)$ by studying the families of surfaces mentioned above.

According to problem (\ref{2}), in \cite{A1} H. Al-Zoubi studied finite type tubes corresponding to the second fundamental form and he proved that: All tubes in $\mathbb{E}^{3}$ are of infinite type. However, for all other classical families of surfaces, the classification of its finite type surfaces is not known yet.

Concerning problem (\ref{3}), ruled surfaces and tubes are the only families studied according to its finite type classification. More specifically, in \cite{A3} authors have shown that all tubes in $\mathbb{E}^{3}$ are of infinite type, while in \cite{A2}, H. Al-Zoubi and others proved that: Helicoids are the only ruled surfaces of finite $III$-type in the 3-dimensional Euclidean space.

In this paper we will pay attention to surfaces of finite $II$-type. First, we will establish a formula for $\Delta^{II}\boldsymbol{x}$ and $\Delta^{II}\boldsymbol{n}$ by using tensors calculations. Further, we continue our study by proving finite type surfaces for an important class of surfaces, namely, ruled surfaces in the Euclidean 3-space.

\section{Preliminaries}

Let $S$ be a (connected) surface in a Euclidean 3-space $E^{3}$ referred to any system of coordinates $u^{1},\ u^{2}$, which does not contain parabolic points, We denote by
\begin{equation*}
I = g_{ij}\,d u^i d u^j, \quad II = b_{ij}\,d u^i d u^j, \quad III =
e_{ij}\,d u^i d u^j, \quad i,j = 1,2,
\end{equation*}
the first, second and third fundamental forms of $S$ respectively. For two sufficiently differentiable functions $f(u^{1}, u^{2})$ and $h(u^{1}, u^{2})$ on $S$, the first differential parameter of Beltrami corresponding to the fundamental form $J = I, II, III$ is defined by \cite{H1}
\begin{equation}  \label{4.1}
\nabla^{J}(f,h):=a^{ij}f_{/i}h_{/j}
\end{equation}
where $f_{/i}: =\frac{\partial f}{\partial u^{i}}$, and $(a^{ij})$ denotes the components of the inverse tensor of $(g_{ij}), (b_{ij})$ and $(e_{ij})$ for $J = I, II$ and $III$ respectively.
The second differential parameter of Beltrami corresponding to the fundamental form $J = I, II, III$ of $S$ is defined by \cite{H1}
\begin{equation}  \label{4.2}
\triangle^{J}f:=-a^{ij}\nabla^{J}_{i}f_{/j},
\end{equation}
where $\nabla^{J}_{i}$ is the covariant derivative in the $u^{i}$ direction corresponding to the fundamental form $J$ and $(a^{ij})$ stands, as in definition (\ref{4.1}), for the inverse tensor of $(g_{ij}), (b_{ij})$ and $(e_{ij})$ for $J = I, II$ and $III$ respectively.

Firstly, we mention the following two relations for later use \cite{S1}:
\begin{equation}  \label{4.21}
\nabla^{II}(h,\boldsymbol{n})+grad^{I}h=0,
\end{equation}
\begin{equation}  \label{4.22}
\nabla^{II}(h,\boldsymbol{x})+grad^{III}h=0,
\end{equation}

Applying (\ref{4.2}) for the position vector $\boldsymbol{x}$ of $S$ we have
\begin{equation}  \label{4.3}
\triangle^{II}\boldsymbol{x}=-b^{ij}\nabla^{II}_{j}\boldsymbol{x}_{/i}.
\end{equation}

Recalling the equations
\begin{equation*}
\nabla^{II}_{j}\boldsymbol{x}_{/i}=-\frac{1}{2}b^{kr}(\nabla^{I}_{k}b_{ij})\boldsymbol{x}_{/r} + b_{ij}\boldsymbol{n},
\end{equation*}
(see \cite{H1}, p.128) and inserting these into (\ref{4.3}), one finds
\begin{equation}  \label{4.4}
\triangle^{II}\boldsymbol{x}=\frac{1}{2}b^{kr}b^{ij}(\nabla^{I}_{k}b_{ij})\boldsymbol{x}_{/r}- b^{ij}b_{ij}\boldsymbol{n},
\end{equation}

From the Mainardi-Codazzi equations (see \cite{H1}, p.128)
\begin{equation}  \label{4.41}
\nabla^{I}_{k}b_{ij}-\nabla^{I}_{i}b_{jk}=0,
\end{equation}
relation (\ref{4.4}) becomes
\begin{equation}  \label{4.5}
\triangle^{II}\boldsymbol{x}=\frac{1}{2}b^{kr}b^{ij}\nabla^{I}_{i}b_{jk}\boldsymbol{x}_{/r}- 2\boldsymbol{n}.
\end{equation}

We consider the Christoffel symbols of the second kind corresponding to the first, second and third fundamental form, respectively
\begin{equation*}
\Gamma^{k}_{ij}:=\frac{1}{2}g^{kr}(-g_{ij/r} + g_{ir/j} + g_{jr/i}),
\end{equation*}
\begin{equation*}
\Pi^{k}_{ij}:=\frac{1}{2}b^{kr}(-b_{ij/r} + b_{ir/j} + b_{jr/i}),
\end{equation*}
\begin{equation*}
\Lambda^{k}_{ij}:=\frac{1}{2}e^{kr}(-e_{ij/r} + e_{ir/j} + e_{jr/i}),
\end{equation*}

and we put
\begin{equation}  \label{4.6}
T_{ij}^{k}: = \Gamma^{k}_{ij} -  \Pi^{k}_{ij},
\end{equation}
\begin{equation}  \label{4.61}
\widetilde{T}_{ij}^{k}: = \Lambda^{k}_{ij} -  \Pi^{k}_{ij}.
\end{equation}

It is known that (see \cite{H1}, p.22)
\begin{equation}  \label{4.7}
T_{ij}^{k}: = -\frac{1}{2}b^{kr}\nabla^{I}_{r}b_{ij},
\end{equation}
\begin{equation}  \label{4.71}
\widetilde{T}_{ij}^{k}: = -\frac{1}{2}b^{kr}\nabla^{III}_{r}b_{ij},
\end{equation}
and
\begin{equation}  \label{4.72}
\widetilde{T}_{ij}^{k}+T_{ij}^{k} = 0.
\end{equation}

Using (\ref{4.6}) and (\ref{4.7}), relation (\ref{4.5}) becomes
\begin{equation}  \label{4.88}
\triangle^{II}\boldsymbol{x}=-b^{kr}T^{j}_{kj}\boldsymbol{x}_{/r}- 2\boldsymbol{n} = -b^{kr}(\Gamma^{j}_{kj}-\Pi^{j}_{kj})\boldsymbol{x}_{/r}- 2\boldsymbol{n}.
\end{equation}

For the Christoffel symbols $\Gamma^{j}_{kj}$ and $\Pi^{j}_{kj}$ we have (see \cite{H1}, p.125)
\begin{equation}  \label{4.8}
\Gamma^{j}_{ij}: = \frac{g_{/i}}{2g}, \ \ \ \ \ \,\Pi^{j}_{ij}:=\frac{b_{/i}}{2b},
\end{equation}
where $g: = det(g_{ij})$ and $b: = det(b_{ij})$. Thus, relation (\ref{4.88}) becomes
\begin{equation}  \label{4.89}
\triangle^{II}\boldsymbol{x}= -\frac{1}{2} b^{kr}(\frac{g_{/k}}{g}-\frac{b_{/k}}{b}) \boldsymbol{x}_{/r}- 2\boldsymbol{n}.
\end{equation}

On the other hand, the Gauss curvature $K$ of $S$ is given by
\begin{equation*}
K=\frac{b}{g}.
\end{equation*}

Once, we have
\begin{equation}  \label{4.81}
\frac{K_{/k}}{K}=\frac{b_{/k}}{b}-\frac{g_{/k}}{g},
\end{equation}
it follows that
\begin{equation}  \label{4.91}
\triangle^{II}\boldsymbol{x}= \frac{1}{2K} b^{kr}K_{/k} \boldsymbol{x}_{/r}- 2\boldsymbol{n}= \frac{1}{2K}\nabla^{II}(K,\boldsymbol{x})- 2\boldsymbol{n}.
\end{equation}

Hence, we obtain, in view of (\ref{4.22}), the following relation
\begin{equation}  \label{4.9}
\triangle^{II}\boldsymbol{x}=-\frac{1}{2K}grad^{III}(K)- 2\boldsymbol{n}.
\end{equation}

We compute now $\triangle^{II}\boldsymbol{n}$.
Taking into consideration the equations (\cite{H1}, p.128)
\begin{equation*}
\nabla^{II}_{i}\boldsymbol{n}_{/j} =-\frac{1}{2}b^{kr}(\nabla^{III}_{r}b_{ij})\boldsymbol{n}_{/k}-e_{ij}\boldsymbol{n},
\end{equation*}
so that
\begin{equation*}
\triangle^{II}\boldsymbol{n}=-b^{ij}(\nabla^{II}_{i})\boldsymbol{n}_{/j},
\end{equation*}
takes the form
\begin{equation*}
\triangle^{II}\boldsymbol{n}=\frac{1}{2}b^{kr}b^{ij}(\nabla^{III}_{r}b_{ij})\boldsymbol{n}_{/k}+ b^{ij}e_{ij}\boldsymbol{n}.
\end{equation*}
On account of
\begin{equation*}
2H = b_{ij} g^{ij} = e_{ij} b^{ij},
\end{equation*}
and (\ref{4.71}) we obtain
\begin{equation*}
\triangle^{II}\boldsymbol{n}=-b^{ij}\widetilde{T}_{ij}^{k}\boldsymbol{n}_{/k}+ 2H\boldsymbol{n}.
\end{equation*}
		
On use of (\ref{4.41}), (\ref{4.7}) and (\ref{4.72}) we have
\begin{equation}  \label{4.10}
\triangle^{II}\boldsymbol{n}=b^{kr}T_{rj}^{j}\boldsymbol{n}_{/k}+ 2H\boldsymbol{n}.
\end{equation}

On the other hand using (\ref{4.6}), (\ref{4.7}), (\ref{4.8}) and (\ref{4.81}) we have
\begin{equation*}
b^{kr} T_{rj}^{j} \boldsymbol{n}_{/k}= -\frac{1}{2K}b^{kr} K_{/r} \boldsymbol{n}_{/k}=-\frac{1}{2K}\nabla^{II}(K,\boldsymbol{n}).
\end{equation*}
		
Inserting this in (\ref{4.10}) we find in view of (\ref{4.21})
\begin{equation}  \label{4.11}
\triangle^{II}\boldsymbol{n}=\frac{1}{2K}grad^{I}(K)+ 2H\boldsymbol{n}.
\end{equation}

From (\ref{4.9}) and (\ref{4.11}) we obtain the following two results which were proved in \cite{S1}:
\begin{theorem}\label{T1}
A surface $S$ in $\mathbb{E}^{3}$ is of finite $II$-type 1 if and only if $S$ is part of a sphere.
\end{theorem}
\begin{theorem}\label{T2}
The Gauss map of a surface $S$ in $\mathbb{E}^{3}$ is of finite $II$-type 1 if and only if $S$ is part of a sphere.
\end{theorem}

Up to now, the only known surfaces of finite $II$-type in $\mathbb{E}^{3}$ are parts of spheres. In the next section we focus our attention on the class of ruled surfaces. Our main result is the following

\begin{theorem}\label{T3}
All ruled surfaces in the three-dimensional Euclidean space are of infinite $II$-type.
\end{theorem}

\section{Proof of Theorem \ref{T3}}

In the Euclidean 3-space $\mathbb{E}^{3}$ let $S$ be a ruled $C^{r}$-surface, $r\geq3$, of nonvanishing Gaussian curvature defined by an injective $C^{r} $-immersion $\boldsymbol{x}=\boldsymbol{x}(s,t)$ on a region $U: = L \times \mathbb{R} \,\,\,(L\subset \mathbb{R}$ open interval) of $\mathbb{R}^{2}$.\footnote {The reader is referred to \cite{P1} for definitions and formulae on ruled surfaces.}
The surface $S$ can be expressed in terms of a directrix curve $\varGamma \colon \boldsymbol{\gamma } = \boldsymbol{\gamma }(s)$ and a unit vector field $\boldsymbol{\rho }(s)$ pointing along the rulings as follows
\begin{equation}  \label{2.1}
S \colon \boldsymbol{x}(s,t) = \boldsymbol{\gamma }(s) + t \, \boldsymbol{\rho }(s), \quad s\in I, t \in \mathbb{R}.
\end{equation}
We further suppose that $\boldsymbol{\rho}$ has arc-length parametrization. Then we have
\begin{equation*}
\langle \boldsymbol{\gamma}^{\prime },\boldsymbol{\rho }\rangle = 0, \quad
\langle \boldsymbol{\rho },\boldsymbol{\rho }\rangle = 1,\quad \langle
\boldsymbol{\rho }^{\prime },\boldsymbol{\rho }^{\prime }\rangle = 1,
\end{equation*}
where the differentiation with respect to $s $ is denoted by a prime and $\langle \,,\rangle$ denotes the standard scalar product in $\mathbb{E}^{3}$.
It is easily verified that the first and the second fundamental forms of $S$ are given by
\begin{align*}
I &= n \, ds^{2} + dt^{2}, \\
II &= \frac{m}{\sqrt{n}}\,ds^{2} + \frac{2A}{\sqrt{n}}\,ds\,dt,
\end{align*}%
where
\begin{align*}
n &= \langle \boldsymbol{\gamma }^{\prime },\boldsymbol{\gamma }^{\prime}\rangle +2\langle \boldsymbol{\gamma }^{\prime },\boldsymbol{\rho }
^{\prime }\rangle t+t^{2}, \\
m &= \left( \boldsymbol{\gamma }^{\prime },\boldsymbol{\rho },\boldsymbol{\gamma }^{\prime \prime }\right) +\left[ \left( \boldsymbol{\gamma }%
^{\prime},\boldsymbol{\rho },\boldsymbol{\rho }^{\prime \prime }\right)
+\left(\boldsymbol{\rho }^{\prime },\boldsymbol{\rho },\boldsymbol{\gamma }^{\prime \prime }\right) \right] t+\left( \boldsymbol{\rho }^{\prime },%
\boldsymbol{\rho },\boldsymbol{\rho }^{\prime \prime }\right) t^{2}, \\
A &= \left( \boldsymbol{\gamma }^{\prime },\boldsymbol{\rho },\boldsymbol{\rho }^{\prime }\right).
\end{align*}
If, for simplicity, we put
\begin{align*}
\zeta &:= \langle \boldsymbol{\gamma }^{\prime },\boldsymbol{\gamma }^{\prime }\rangle ,\qquad \eta :=\langle \boldsymbol{\gamma } ^{\prime },%
\boldsymbol{\rho }^{\prime }\rangle , \\
\mu &:=\left( \boldsymbol{\rho }^{\prime },\boldsymbol{\rho },\boldsymbol{\rho }^{\prime \prime }\right),\quad \nu :=\left( \boldsymbol{\gamma }%
^{\prime },\boldsymbol{\rho },\boldsymbol{\rho }^{\prime \prime}\right)+\left( \boldsymbol{\rho }^{\prime },\boldsymbol{\rho },\boldsymbol{%
\gamma}^{\prime \prime }\right),\quad \xi :=\left( \boldsymbol{\gamma }^{\prime },\boldsymbol{\rho },\boldsymbol{\gamma }^{\prime \prime }\right),
\end{align*}%
we have
\begin{equation*}
n = t^{2} + 2\eta \, t + \zeta , \quad m = \mu \, t^{2} + \nu \,t + \xi .
\end{equation*}

For the Gauss curvature $K$ of $S$ we find
\begin{equation*}
K = -\frac{A^{2}}{n^{2}}.
\end{equation*}
The second Beltrami differential operator with respect to the second fundamental form after a long computation is given by \cite{Y1}
\begin{equation}\label{2.2}
\triangle^{II} =- \frac{\sqrt{n}}{A}\bigg( -2\frac{\partial^{2}}{\partial s\partial t}+\frac{m}{A}\frac{\partial^{2}}{\partial t^{2}}+\frac{m_{t}}{A}\frac{\partial}{\partial t}\bigg),
\end{equation}
where $m_{t}:=\frac{\partial m}{\partial t}$.

Applying (\ref{2.2}) for the position vector $\boldsymbol{x}$, it follows:
\begin{eqnarray}\label{3.2}
\triangle ^{II}\boldsymbol{x} & = &-\frac{1}{\sqrt{n}}\bigg(-\frac{2n}{A}\boldsymbol{\rho }^{\prime}+\frac{nm_{t}}{A^{2}}\boldsymbol{\rho }\bigg)=\frac{1}{\sqrt{n}}\boldsymbol{P_{1}(t)}
\end{eqnarray}%
where $\boldsymbol{P_{1}(t)}$ is a vector-valued function in $\mathbb{E}^{3}$ whose components are polynomials in $t$ of degree less than or equal 3 with functions in $s$ as coefficients. More precisely, we have
\begin{equation*}
\boldsymbol{P_{1}(t)}=\frac{1}{A^{2}} \big[2\mu\boldsymbol{\rho }t^{3}+\big((4\mu\eta+\nu)\boldsymbol{\rho }+2A\boldsymbol{\rho }^{\prime}\big)t^{2}  \\
+\big((2\zeta\mu+2\eta\nu)\boldsymbol{\rho }+4\eta A\boldsymbol{\rho }^{\prime}\big)+(\zeta\eta\boldsymbol{\rho }+2\zeta A\boldsymbol{\rho }^{\prime})\big].
\end{equation*}

Before we start the proof of our theorem we give the following Lemma which can be proved after a somewhat long but straightforward calculation.
\begin{lemma}\label{L2.1}
Let $g$ be a polynomial in $t$ of degree $d$ with functions in $s$ as coefficients. Then $\triangle ^{II} \big(\frac{g}{n^{r}}\big) =\frac{\widehat{g}}{n^{r+\frac{3}{2}}}$,
where $\widehat{g}$ is a polynomial in $t$ with functions in $s$ as coefficients and $\deg(\widehat{g})\leq d + 4$.
\end{lemma}

If the ruled surface $S$ is of finite $II-type$, then for some natural number $k$ there exist real numbers $c_{1}, c_{2}, \dotsb  , c_{k}$ such that
\begin{equation}  \label{e5}
\left(\triangle ^{II}\right)^{k+1}\,\boldsymbol{x} + c_{1}\left(\triangle^{II}\right)^{k} \, \boldsymbol{x} + \dotsb + c_{k}\, \triangle ^{II}\,%
\boldsymbol{x} = \mathbf{0},
\end{equation}%
see \cite{A3}. By applying Lemma \ref{L2.1}, we conclude that there is an $\mathbb{E}^{3}$-vector-valued function $\boldsymbol{P}_{k}$ in the variable $t$ with some functions in $s$ as coefficients, such that

\begin{equation*}
\left(\triangle ^{II}\right)^{k} \, \boldsymbol{x} = \boldsymbol{P}_{k}(t),
\end{equation*}%
where $\deg(\boldsymbol{P}_{k})\leq 4k-1$ and $r = \frac{3}{2}k-1$. Now, if $k$ goes up by one, the degree of each component of $\boldsymbol{P}_{k}$ goes up at most by 4 while the degree of the denominator goes up by  $\frac{3}{2}k-1$. Therefore, the sum (\ref{e5}) can never be zero, unless
\begin{equation}  \label{tr}
\triangle ^{II}\,\boldsymbol{x} = \boldsymbol{P}_{1}=\mathbf{0}.
\end{equation}

But then
\begin{equation}  \label{2.5}
-2\boldsymbol{\rho }^{\prime}+\frac{m_{t}}{A}\boldsymbol{\rho }=\mathbf{0}.
\end{equation}

By taking the derivative of $\langle \boldsymbol{\rho },\boldsymbol{\rho }\rangle = 1,$ we observe that the vectors $\boldsymbol{\rho }$ and $\boldsymbol{\rho }^{\prime}$ are linearly independent. Thus (\ref{2.5}) cannot be achieved unless $\boldsymbol{\rho }$ is constant, which implies that $K$ $\equiv 0$. This is clearly impossible for the surfaces under consideration. The proof of the theorem is completed.

\end{document}